\documentclass{mymathart}
\usepackage
           {mymathmacros}

\usepackage{mathdots}

\newtheorem{thm}{Theorem}[section]

\newtheorem{cor}[thm]{Corollary}
\newtheorem{prop}[thm]{Proposition}
\newtheorem{rmk}[thm]{Remark}

\newcommand{\Ek}{E_{\kappa}}

\newcommand{\Sequ}{\mathcal{S}}

\newcommand{\adds}{\underline{s}}

\newcommand{\s}{\adds}

\newcommand{\ts}{t_{\s}}

\newcommand{\addst}{\underline{\tilde{s}}}
\newcommand{\st}{\widetilde{s}}

\newcommand{\gs}{g_{\adds}}
\newcommand{\gsk}{\gs^{\kappa}}

\newcommand{\I}{\mathcal{I}}

\newcommand{\PR}{G_{\s}}
\newcommand{\PRS}[1]{G_{#1}}

\newcommand{\PRconj}{\mathcal{G}}
\newcommand{\DRconj}{\mathfrak{g}}

\newcommand{\M}{\mathcal{M}}

\newcommand{\Ih}{\widehat{\I}}

\newcommand{\sm}{\setminus}

\usepackage{hyperref}

\title{Classification of Escaping Exponential Maps}

\author{Markus F\"orster}
\address{
KPMG Deutsche Treuhand-Gesellschaft,
Marie-Curie-Stra{\ss}e 30,
60439 Frankfurt/Main, Germany}
\email{mfoerster@kpmg.com}
\author{Lasse Rempe}
\address{Dept.~of Mathematical Sciences, University of Liverpool, 
  Liverpool L69 7ZL,
United Kingdom}
\email{l.rempe@liverpool.ac.uk}
\thanks{The first author was supported in part
 by a European fellowship of the Marie Curie Fellowship Association.
 The second author was supported in part
 by a postdoctoral fellowship of the 
 German Academic Exchange Service (DAAD) and
 by the German-Israeli Foundation
 for Scientific Research and Development (G.I.F.),
 grant no.\ G-643-117.6/1999}
\author{Dierk Schleicher}
\address{International University Bremen, P.O.~Box 750 561, 
28725 Bremen, Germany}
\email{dierk@iu-bremen.de}


\date{\today}

\begin{document}

\begin{abstract}
 We give a complete classification of the set of parameters $\kappa$ for
 which the singular value of
 $E_{\kappa}:z\mapsto \exp(z)+\kappa$ escapes to $\infty$ under
 iteration. In particular, we show that every path-connected component
 of this set is a curve to infinity.
\end{abstract}

\maketitle

\section{Introduction}

In 1777, Euler \cite{euler} 
 studied the question for which $a > 0$ the limit
 $a^{a^{a^{\iddots}}}$
 exists. Writing $\lambda=\ln a$, this amounts to asking for which
 $\lambda\in\R$ the iteration 
 \begin{equation}
   \label{eqn:euler}
   \lambda_0:=0, \quad \lambda_{n+1}:=e^{\lambda \cdot \lambda_n}
 \end{equation}
  has a limit, which 
  is easily seen to be the case if and
  only if $-e \leq \lambda \leq 1/e$.
  If we instead allow $\lambda$ to vary in the complex plane, 
  the question obtains a much richer structure. 
  It is readily verified that
  $(\lambda_n)$ has a nontrivial limit if and only if 
  $\lambda$ is of the form $\lambda=\mu e^{-\mu}$, where
  $|\mu|<1$ or $\mu$ is a root of unity. In addition, 
  there are countably many 
  values of $\lambda$
  for which the sequence is eventually constant; a complete
  classification of these can be found in 
  \cite{misiurewiczclassification,expmisiurewicz}.

 This exhausts all cases for which there is a limit in $\C$, but
  when does $\lambda_n\to\infty$?
  This is certainly the case for all real
  $\lambda>1/e$. However, the complete answer is much more interesting:
 \begin{thm} \label{thm:main}
  The set of parameters $\lambda\in\C$ for which $\lambda_n\to\infty$
  consists of
  uncountably many disjoint curves in $\C$. 
  More precisely, every path-connected component of this set
  is an injective curve $\gamma:(0,\infty)\to\C$ or
  $\gamma:[0,\infty)\to\C$ with $\lim_{t\to\infty} \gamma(t) = \infty$.
 \end{thm}

 This result is of interest not only in its own right, 
  but provides an important ingredient in the dynamical study 
  of the family of complex exponential maps. In this article, 
  we will parametrize this family as
  \[ E_{\kappa}:\C\to\C; z\mapsto \exp(z)+\kappa \] 
  rather than $\exp(\lambda z)$ as above\footnote{%
  We prefer our parametrization since it affords all
  maps the same asymptotic behavior near infinity.}%
 \ (note that the two maps
   are conformally conjugate to each other 
   if $\lambda = \exp(\kappa)$). 
  Theorem \ref{thm:main} will thus be proved by
  showing the corresponding statement for the set
  \[ \I := \{\kappa\in\C: \Ek^n(\kappa)\to\infty \} \]
  of \emph{escaping parameters}.
 
 This set is an analog to the exterior of the
  Mandelbrot set $\M$ for quadratic polynomials. The foliation of
  the open set
  $\C\setminus\M$ by \emph{parameter rays} introduced by
  Douady and Hubbard
  \cite{orsay} has played a vital role in the understanding of 
  quadratic polynomials. Although $\I$ is not
  an open set (in fact, it is nowhere dense), it has long
  been known to contain curves which can play the role
  of parameter rays (compare the historical remarks below).
  We will in fact
  provide a complete classification of escaping parameters
  in terms of these rays (Corollary \ref{cor:classification}):
  \emph{every escaping parameter is either on a parameter ray or the
  endpoint of such a ray}. 
  This is analogous to (and builds on) a similar classification 
  \cite{expescaping} obtained
  for the \emph{set of escaping points},
  \[ I(\Ek) := \{z\in\C: \Ek^n(z)\to\infty \text{ as $n\to\infty$}\}, \]
  of an arbitrary exponential map $\Ek$, in terms of \emph{dynamic rays}.

 Our description of the set $\I$
  also shows that the 
  ``Dimension Paradox'' discovered in \cite{karpinska} 
  occurs in the parameter plane as well. 
  Indeed, Qiu \cite{hausdorffparameter} showed that
  $\I$ has Hausdorff dimension
  $2$. On the other hand, the
  set of parameter rays (without endpoints)
  only has Hausdorff dimension $1$
  \cite{parameterdimensionparadox2}. Thus $\I$ contains a
  two-dimensional set of endpoints, each connected to
  $\infty$ by a curve in $\I$, such that all curves are
  disjoint and the union of these curves
  has Hausdorff dimension one. 

\smallskip

\noindent
\textsc{Historical remarks. }
Parameter space of exponential maps was first studied in
 \cite{bakerexp,alexmisharussian}, though the focus of these articles
 lay on the structure of hyperbolic components. A
 construction of 
 (tails of) certain parameter rays
 consisting of escaping parameters was first sketched in
 \cite{devaneybifurcation}, and carried out in
 \cite{dgh}. (Compare the remark after Proposition
  \ref{prop:markus}.)
  In \cite{habil}, full parameter rays
  were constructed 
  for arbitrary bounded external addresses
  and in \cite{markus,markusdierk}, this was generalized
  to all parameter rays. 

 Classically, both dynamical and parameter
  rays were often referred to as ``hairs'' in the transcendental setting.
  We chose our terminology to emphasize the analogy with the polynomial
  case. Similarly, we chose the term \emph{external address},
  introduced in \cite{expescaping}, over that of \emph{itinerary}.
  On the one hand, this highlights the relation to external angles
  of polynomials; on the other, it avoids the overloading of the
  latter term, which is frequently used for another kind of
  symbolic sequence
  both in exponential and polynomial dynamics
  (compare e.g.\ \cite[Section 4.1]{expper}).

\section{Preliminaries} \label{sec:preliminaries}

 For the remainder of this paper, we will fix the function 
  \[ F:[0,\infty)\to [0,\infty);t\mapsto \exp(t)-1 \] 
  as a model for exponential growth. 
  Following \cite{expescaping}, a sequence $\s= s_1 s_2 s_3 \dots$ of 
  integers is called an \emph{external address}; we also define 
  the \emph{shift map} on external addresses
  by $\sigma(s_1 s_2 s_3 \dots ) := s_2 s_3 \dots$. 

  If $\Ek$ is an 
  exponential map and $z\in \C$, then we say that $z$ \emph{has external 
  address $\s$} if $\im \Ek^{k-1}(z) \in 
  \bigl((2s_{k}-1)\pi,(2s_{k}+1)\pi\bigr)$ for all $k\geq 1$. An 
  external address $\s$ is called \emph{exponentially bounded} if there is 
  an $x>0$ such that 
  $2\pi|s_k| < F^{k-1}(x)$ 
  for all $k\geq 1$; the set of all such addresses 
  is denoted by $\Sequ_0$.
  An address $\s\in \Sequ_0$ is called \emph{slow} if there exists
  a subsequence $\sigma^{n_k}(\s)$ of the $\sigma$-orbit of $\s$
  such that
  all $\sigma^{n_k}(\s)$
  are exponentially bounded
  with the same $x$.
  Otherwise, $\s$ is called \emph{fast}. 
  For any external address $\s$, let us also define
   $\ts := \limsup_{k\to\infty} F^{-(k-1)}(2\pi|s_{k}|)$, where
   $F^{-(k-1)}$ denotes the $(k-1)$-th iterate of the function 
   $F^{-1}:[0,\infty)\to [0,\infty)$. Note that
   $\Sequ_0 = \{\s\in\Z^{\N}: \ts < \infty\}$.

\begin{figure}%
 \begin{center}%
  \resizebox{.98\textwidth}{!}{%
  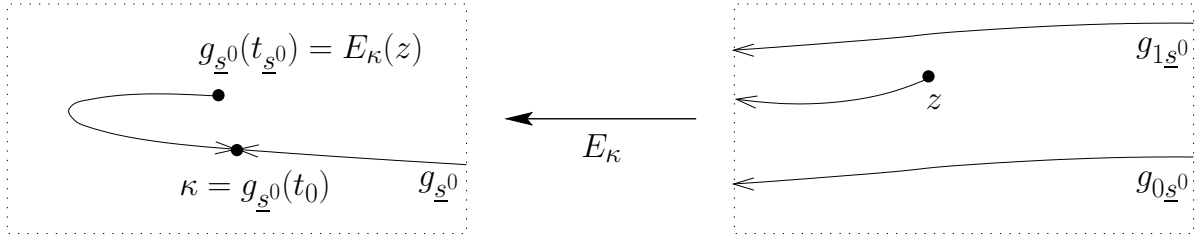}%
 \end{center}%
 \caption{An illustration of the second part of Proposition
   \ref{prop:classification}. Here we have $\kappa\in g_{\s^0}$ for a
   fast address $\s^0$. Hence, its
   preimage rays (at addresses $k\s^0$, $k\in\Z$) end
   prematurely, and the
   preimage $z$ of $g_{\s^0}(t_{\s^0})$ is not
   itself associated to any dynamic ray.\label{fig:brokenray}}
\end{figure}

\begin{prop}[Dynamic rays {\cite[Corollary 6.9]{expescaping}}] \label{prop:classification} 
 Let $\kappa\in\C$ with
  $\kappa\notin I(\Ek)$. 
 For every $\s\in\Sequ_0$, 
   there exists a unique 
  simple curve $\gs=\gsk:X_{\s}\to I(\Ek)$, where
  $X_{\s}:=(\ts,\infty)$ if $\s$ is slow and 
  $X_{\s}:=[\ts,\infty)$ if $\s$ is fast, such that
 \begin{itemize}
  \item $\gsk(t)$ has external address $\s$ for all sufficiently large $t$ and
  \item for all $t>\ts$, we have $\re \Ek^k(\gsk(t)) =  F^k(t)+o(1) $ as
          $k\to\infty$.
 \end{itemize}

For every
  $z\in I(\Ek)$, there exist unique $\s$ and $t$ such that
  $z = \gsk(t)$. Furthermore,  
  $\Ek(\gs(t)) = g_{\sigma(\s)}(F(t))$ for all $\s$ and $t$.

 If $\kappa\in I(\Ek)$, then the above is still true
  with the modification that there are countably many
  addresses $\s$ such that $\gs$ is not defined
  on the complete interval $X_{\s}$. More precisely,
  there exist some $\s^0\in\Sequ_0$ and $t^0$ such
  that
  $\kappa=g_{\s^0}^{\kappa}(t_0)$. If $\s$ is an address with
  $\sigma^n(\s)=\s^0$, where $n\geq 1$, then 
  $g^{\kappa}_{\s}(t)$ is not defined for
  $t\leq F^{-n}(t_0)$. 
  If $z\in I(\Ek)$, then either $z=\gs(t)$ for some $\s\in\Sequ_0$ and
  $t\in X_{\s}$ (in which case $\s$ and $t$ are unique), or
  $\Ek^n(z) = g_{\s^0}(t)$ for some $n\geq 1$ and $t\leq t_0$.
   \qedd
\end{prop}

 The \emph{dynamic ray} at address $\s$ is defined to be the curve 
  $\gsk\colon (t^*,\infty)\to\C$ for the least possible value of $t^*$. 
  (This will be $t^*=\ts$, unless $\kappa\in I(\Ek)$ and $\s$ is one of the 
  countably many exceptional addresses.) If $t^*=\ts$ and $\s$ is fast, then
  $\gsk(\ts)$ is called the \emph{endpoint} of the ray $\gsk$.

  
\begin{remark}[Remark 1]
 See Figure \ref{fig:brokenray} for an illustration of the last part of
  the Proposition, in the case where
  $\s^0$ is fast. The case where $\s^0$ is slow, but not periodic,
  is similar (except that the ray may not have a landing point 
  in this case). If $\s^0$ is periodic, say of period $n$,
  then $g_{\s^0}$ is itself defined only for $t\geq F^{-n}(t_0)$,
  and every iterated 
  preimage component
  of $g_{\s^0}\bigl((F^{-n}(t_0),t_0)\bigr)$ is an arc
  tending to infinity in both directions (compare
  \cite[Figure 2]{topescapingnew}).
\end{remark}
\begin{remark}[Remark 2]
 We say that the ray $\gs:(t^*,\infty)\to\C$ \emph{lands} at
  $z\in\Ch$ if $z = \lim_{t\to t^*}(\gs(t))$. 
 Proposition \ref{prop:classification} shows that every dynamic ray at a 
  fast address lands at an escaping point, but does not make any statements
  about the accumulation behavior of rays at slow addresses. Although it
  is true that many rays at slow addresses also land
  (compare \cite[Proposition 6.11]{expescaping} and 
  \cite[Theorem 1]{landing2new}),
  this problem is much more complicated in general. Even for
  postsingularly finite exponential maps, there are dynamic rays at
  bounded addresses which do \emph{not} land; compare
  \cite{nonlanding}. 
\end{remark}
\begin{remark}[Remark 3]
 Let us define a topological space 
  $X:=\bigcup_{\s\in\Sequ_0}\left( \{\s\}\times X_{\s}\right)$ 
  (that is, $X$ is the disjoint union of all $X_{\s}$). 
  $X$ is endowed with the product topology of $\Sequ_0\times [0,\infty)$,
  where $\Sequ_0$ is understood as a discrete space and $[0,\infty)$ comes
  with the usual topology.
  Then the first
  part of the proposition implies that for every $\kappa\in\C\sm \I$, 
  there is a
  continuous bijection $\DRconj^\kappa\colon X\to I(E_\kappa)$ given by
  $\DRconj^\kappa(\s,t)=g_{\s}^\kappa(t)$. 

 The topology on $X$ is much finer than the planar topology on
  $I(E_{\kappa})$, and it is impossible to embed
  $X$ in the plane.  Thus $\DRconj$ can never be a homeomorphism. 
  However, it is possible 
  \cite{topescapingnew} to endow $X$ with a different topology%
   \footnote{It is useful to think of this new topology 
    as induced by
    the product topology of the lexicographic order on
    $\Sequ_0$ and the usual topology on the real line. However, this is
    not quite accurate with our definitions: while it yields the
    correct topology in the $\s$-direction, we would need to change the
    parametrization of dynamic rays to obtain the right notion of
    convergence of potentials. See \cite{topescapingnew} for details.}
\ ---
  turning it into a \emph{straight brush} in the sense of
  \cite{aartsoversteegen} --- such that $\DRconj$ is a homeomorphism when
  restricted to large subsets of $X$ (compare the proof of Corollary
  \ref{cor:dynamicrays}).
  When the singular value $\kappa$ belongs to the Fatou set, then
  in fact $I(\Ek)$ is a straight brush \cite{tying}, and the map
  $\DRconj$ is a 
  homeomorphism between the straight brush $X$ and $I(\Ek)$
  \cite{topescapingnew}. It seems
  reasonable to expect that $I(\Ek)$ is \emph{not} a straight brush
  whenever
  $\kappa$ belongs to the Julia set; indeed, for a large class of
  such parameters it is known that there are dynamic rays 
  which do not land and whose accumulation
  sets contain escaping points \cite{nonlanding}. 
\end{remark}

Let us now turn to the question how escaping parameters 
 are organized in the parameter plane. By the above
 result, for every such parameter $\kappa$ there exist unique
 $\s$ and $t$ such that 
 $\kappa=\gsk(t)$. The set of parameters for which 
 $t>\ts$ --- i.e.\ for which the singular value is not the endpoint of
 a dynamic ray --- 
 is described by the following result.

\begin{prop}[Parameter rays \cite{markus,markusdierk}]  \label{prop:markus}
 Let $\s\in\Sequ_0$.
 Then for every
 $t>\ts$, there exists a unique parameter $\kappa=:\PR(t)$ with
 $\kappa= \gsk(t)$. Furthermore,
 $\PR:(\ts,\infty)\to\C$ is continuous, and
  $|\PR(t)\ - (t+ 2\pi i s_1)|\to 0$
 as $t\to\infty$. \qedd
\end{prop}

The curves $\PR\colon(\ts,\infty)\to\C$ are called
 \emph{parameter rays}. As usual, the ray $\PR$ is said to
 \emph{land} at a point $\kappa\in\C$ if $\PR(t)\to\kappa$ as $t\searrow\ts$.

\begin{remark}[Remark 1]
  Proposition \ref{prop:markus} was proved in \cite{dgh} for the
 case of those bounded external addresses which are ``regular''.
 (With our notation, these are exactly those addresses satisfying
 $s_k\neq s_1$ for $k>1$.) Also, for \emph{every} bounded address $\s$, a
 \emph{tail} of the parameter ray $\PR$ was constructed; that is, the
 existence of a curve $\PR\colon(\ts',\infty)\to\C$ was proved for sufficiently large $\ts'\ge\ts$.
\end{remark}
\begin{remark}[Remark 2]
 For parameter rays
  at regular bounded addresses, \cite{dgh} also proves the 
  existence of (non-escaping)
  landing points. In the present article, we are studying
  only the structure of the set $\I$ of escaping parameters,
  and are hence not making any claims regarding the landing
  of parameter rays at general slow addresses, which is 
  a much more difficult problem. Indeed, 
  one can ask whether it is true that every parameter ray
  lands, and that every point in $\cl{\I}$ (the \emph{bifurcation locus})
  is either on a ray or a landing point of a ray. This is
  an analog of the celebrated and unresolved
  \emph{MLC conjecture} for quadratic
  polynomials, and would imply density of hyperbolicity in exponential
  parameter space by \cite[Corollary A.5]{expcombinatorics}. 
\end{remark}

In order to extend Proposition \ref{prop:markus} to escaping endpoints, 
 we shall need
 two ingredients. The first is the fact that escaping endpoints depend
 holomorphically on the parameter.
 \begin{prop}[Holomorphic dependence \protect{\cite[Corollary 4.4]{topescapingnew}}]
  \label{prop:holomorphicdependence}
   Let $\s$ be a fast external address, and suppose that $\kappa$ is a
   parameter with 
   $\kappa\notin g^{\kappa}_{\sigma^n(\s)}([t_{\sigma^n(s)},\infty))$ 
   for all $n\geq
   1$. Then there exists a neighborhood $U$ of $\kappa$ such that the
   function
    \[ U\times [\ts,\infty) \to \C, (\kappa,t)\mapsto \gsk(t) \]
   is defined, continuous in both
   variables and holomorphic in $\kappa$ for fixed $t$. \qedd
  \end{prop}
\begin{remark}
 In \cite{topescapingnew}, a different parametrization of $\gsk$ was
 used. This is not essential because the reparametrization which
 occurs is independent of $\kappa$.
\end{remark}

We also require the following ``Squeezing
 Lemma'' which was proved in \cite{boundary}. 
\begin{prop}[{\cite[Theorem 3.5]{boundary}}]
 Let $\gamma:[0,\infty)\to\C$ be a curve in parameter space
 such that for all $t$,
 all periodic points of $E_{\gamma(t)}$ are repelling. 
 If furthermore
 $\lim_{t\to\infty}|\gamma(t)|=\infty$, then $\gamma$ contains the curve
 $\PR\bigl([t_0,\infty)\bigr)$ for some 
  $\s\in\Sequ_0$
 and $t_0>\ts$. \qedd
\end{prop}

\section{Classification of Escaping Parameters}

We are now ready to give a classification of those parameters
 for which the singular value is the escaping endpoint of
 a dynamic ray.    

\begin{thm} \label{thm:endpointclassification}
 Let $\s\in\Sequ_0$ be fast. Then there is a unique parameter 
  $\kappa$ with $\kappa=\gsk(t_{\s})$. Furthermore, the parameter ray
  $\PR$ lands at this parameter $\kappa$. 
\end{thm}
\begin{proof} 
 Let $L_{\s}$ denote 
  the limit set of
  the parameter ray $\PR$ as $t\searrow t_{\s}$, and
  let $A_{\s}$ consist of all parameters $\kappa\in\C$ with
  $\gs^{\kappa}(\ts)=\kappa$.
  We need to show that $L_{\s}=A_{\s}$ and
  that this set consists of a single parameter.

 In view of the results cited in the previous section, the idea of
  the proof is relatively straightforward. By continuity of
  $\gs^{\kappa}(t)$ in $\kappa$ and $t$, we ought to have $L_{\s}\subset
  A_{\s}$; as the zero set of an analytic function,
  $A_{\s}$ is discrete. Hence $L_{\s}$ consists of a single point,
  and $A_{\s}\subset L_{\s}$ follows from Hurwitz's theorem. However, 
  a little bit of care is required to carry this idea through, as
  e.g.\ it is a priori
  conceivable that $L_{\s}$ contains parameters for which
  $\gs^{\kappa}(\ts)$ is not defined.

\begin{claim}[Claim 1]
  $A_{\s}$ is closed and contains
   $L_{\s}$.
\end{claim}
\begin{subproof}
   Consider the set
    \[
       B_{\s} := \bigcap_{t_0>\ts} \cl{\{\kappa\in\C:
                          \exists t\in [\ts,t_0]: \kappa=\gs^{\kappa}(t)\}}
           = \cl{A_{\s}}\cup 
               \bigcap_{t_0>\ts}
                \cl{\PR\bigl((\ts,t_0]\bigr)}
               = \cl{A_{\s}} \cup L_{\s}.
    \]
   To prove Claim 1, we 
    show that
    $B_{\s}\subset A_{\s}$. 
  So let $\kappa\in B_{\s}$ and choose sequences
  $\kappa_j \to \kappa$ and $t_j\to\ts$ with $\gs^{\kappa_j}(t_j)=\kappa_j$.
  Let $n\geq 0$ be the smallest number for
  which
  $g^{\kappa}_{\sigma^n(\s)}(F^n(\ts))$ is defined (such an $n$ exists
  by Proposition \ref{prop:classification}).
  By Proposition \ref{prop:holomorphicdependence},
   \begin{align} \label{eqn:postsingularendpoint}
       g^{\kappa}_{\sigma^n(\s)}(F^n(\ts)) &=
       \lim_{j\to\infty} g^{\kappa_j}_{\sigma^n(\s)}(F^n(t_j)) \\
       &=
       \lim_{j\to\infty} E_{\kappa_j}^n(g^{\kappa_j}_{\s}(t_j))=
       \lim_{j\to\infty} E_{\kappa_j}^n(\kappa_j) =
            E_{\kappa}^n(\kappa).\notag
   \end{align}

 We claim that $n=0$.
  Suppose not. By minimality of
  $n$, it follows that the curve $g^{\kappa}_{\sigma^n(\s)}$ contains
  $\kappa$, say
  $\kappa=g^{\kappa}_{\sigma^n(\s)}(t)$ for some $t \geq  F^n(\ts)$.
  Since $\kappa\neq \Ek^n(\kappa)$, it follows that
   $t>F^n(\ts)$. Now
   \[ g_{\sigma^{2n}(\s)}^{\kappa}(F^n(t)) = E_{\kappa}^n(\kappa) =
      g^{\kappa}_{\sigma^n(\s)}(F^n(\ts))\,\,, \] 
  which contradicts the statement on uniqueness
  of address and potential in
  Proposition
  \ref{prop:classification}. So
  $n=0$, and (\ref{eqn:postsingularendpoint}) reads
  $g^{\kappa}_{\s}(\ts)=\kappa$ as claimed.  
 \end{subproof}

\begin{claim}[Claim 2]
  $A_{\s}$ is discrete. 
\end{claim}
\begin{subproof}
 Suppose, by contradiction, that
  $A_{\s}$ has an accumulation point $z\in \C$; then
  $z\in A_{\s}$ by Claim 1.
  Now $A_{\s}$ is the zero set of
  the function $\kappa\mapsto \gs^{\kappa}(\ts)-\kappa$.
  By Proposition \ref{prop:holomorphicdependence}, 
   this function is analytic and its 
   domain of definition $D$ is open.
   Let $U$ be the connected component of $D$
   containing $z$.
   Then $U\subset A_{\s}\subset D$ by the identity theorem. Thus
   $U$ is a connected
   component both of the open set $D$ and the closed set
   $A_{\s}$, and therefore
   $U=D=A_{\s}=\C$. This is a contradiction as
   $\C\setminus A_{\s}$ is clearly nonempty.
\end{subproof}

  Since the limit set $L_{\s}\subset A_{\s}$ is connected and $A_{\s}$
  is discrete, $L_{\s}$ is either empty or consists of a single
  point. In the former case, the ray $\PR$ would land at
  infinity, which is impossible by the Squeezing Lemma. Thus
  $\PR$ lands at some point in $A_{\s}\subset\C$.

  It remains to show that $A_{\s}\subset L_{\s}$. Let
  $\kappa_0\in A_{\s}$, and let $V$ be any small 
  disk neighborhood
  of $\kappa_0$ with $\cl{V}\cap A_{\s}=\{\kappa_0\}$. If $t>\ts$ is small
  enough, then by Hurwitz's theorem the function $\kappa\mapsto 
  \gs^{\kappa}(t)-\kappa$ has at least one zero in $V$; i.e.\ 
  $\PR(t)\in V$.
\end{proof}

\begin{remark}
 In a similar way, the Squeezing Lemma can 
  be used to give an alternative
  proof of Proposition \ref{prop:markus};
  see \cite[Section 5.12]{thesis}.
\end{remark}

 Together with Propositions \ref{prop:markus} and \ref{prop:classification},
  the previous theorem yields the desired classification of all escaping
  parameters. (Recall that $X_{\s}=[\ts,\infty)$ when $\s$ is a fast external 
  address, and $X_{\s}=(\ts,\infty)$ when $\s$ is slow.)
  
\begin{cor} \label{cor:classification}
 Let $\s\in\Sequ_0$,
 and let $t\in X_{\s}$. Then
 there exists a unique parameter $\kappa=:\PR(t)$ with
 $\kappa=\gsk(t)$. Furthermore, $\PR$ is continuous in $t$; in particular,
 if $\s$ is fast, then $\PR$ lands at an escaping parameter. 

 Conversely, 
  for any escaping parameter $\kappa$, there exist
  unique $\s\in\Sequ_0$ and $t\in X_{\s}$ such that $\kappa=\PR(t)$. \qed
\end{cor}
\begin{remark}
 With the notation introduced after Proposition \ref{prop:classification},
  we can reformulate Corollary \ref{cor:classification} as follows:
  there exists a continuous bijection $\PRconj:X\to \I$ such that,
  for all $x\in X$, the parameter $\kappa := \PRconj(x)$ satisfies
  $\DRconj^{\kappa}(x) = \kappa$. 
\end{remark}

\begin{cor}[First Entry of External Address]
 Let $\kappa_0$ be an escaping parameter, say $\kappa_0=\gs^{\kappa_0}(t)$.
 Then $\im\kappa_0\in \bigl((2s_1-1)\pi,(2s_1+1)\pi\bigr)$.
\end{cor}
\begin{proof} For sufficiently large $t$,
 $\PR(t)$ belongs to the strip
 $\bigl((2s_1 - 1)\pi,(2s_1+1)\pi\bigr)$ by Proposition
 \ref{prop:markus}. Since the lines $\{\im \kappa = (2k+1)\pi\}$ only contain
 attracting and parabolic parameters \cite{bakerexp}, the entire
 parameter ray $\PR$ belongs to this strip. 
\end{proof}

\section{Path Components of Sets Obtained as a Union of Curves}

 To complete the proof of Theorem \ref{thm:main}, we still need to show that
  the sets $G_{\s}(X_{\s})$ indeed form the path-connected components of
  $\I$. To show this fact, we require
  a general topological principle 
  (Proposition \ref{prop:pathcomponentsarecurves}),
  which will be proved
  in this section. For completeness, we also apply this
  principle to show that every dynamic ray (possibly with endpoint) 
  of an exponential map is a path-connected component of the set 
  of escaping points, as
  stated in \cite{expescaping}
  without proof.

 We shall require the following basic fact of continuum theory.

 \begin{prop}
   \label{prop:sigmaconnectedness}
   Let $K$ be a continuum, and let $(A_n)$ be an increasing sequence
    of nonempty closed subsets of $K$ such that
    $K=\bigcup A_n$. 
    If every connected component of
    $A_n$ is also a connected component of $A_{n+1}$, then
    $A_n=K$ for all $n$.
 \end{prop}
 \begin{proof} We will apply this result only when 
  $K=[0,1]$; let us restrict to this case for simplicity. 
  (The proof in the
   general case is similar; see
     \cite[Theorem 5.16]{continuumtheory} for details.)
   If $A_n=K$ for any $n$, then by assumption 
   $A_n=K$ for all $n$. So
   suppose by contradiction that $A_n\neq K$ for all $n$. 
   We
   will inductively
   define a nested
   sequence $(K_j)_{j\geq 0}$ of nonempty compact intervals such that
   $K_j\not\subset A_k$ for all $k$ and such that
   $K_{j+1}\subset K_j\setminus A_{j+1}$.
   This yields a contradiction since the intersection 
   of all $K_j$ is nonempty, but is also disjoint from all $A_j$. 

  We begin the construction by setting $K_0 := K$.   
  To construct $K_{j+1}$ from $K_j$, pick
   some point
   $x\in K_j\setminus A_{j+1}$ (this is possible by
   the induction hypothesis). Then there is some $n > j$ such that
   $x\in A_n$, and the connected component $C$ of $A_n$
   containing $x$ is disjoint from $A_{j+1}$ by hypothesis. Since
   $C$ is compact, we can find a compact interval
   $K_{j+1}$ with 
   $C\subsetneq K_{j+1} \subset K_j\setminus A_{j+1}$. 
   For every $k\in\N$, either $C\not\subset A_k$ or
   $C$ is a connected component of $A_k$. In either case,
   $K_{j+1}\not\subset A_k$, as required. 
 \end{proof}

 \begin{prop} \label{prop:pathcomponentsarecurves}
  Let $I$ be a Hausdorff topological space. Let $\Gamma$ be
   a partition of $I$ into path-connected subsets such that
   no union of two different elements of $I$ is path-connected.

  Suppose that $I$ can be written as a countable union of closed subsets
   $I_k$ such that
  \begin{enumerate}
   \item $I_{k}\subset I_{k+1}$
   \item every path-connected component of $I_k$ is
     contained in some element of $\Gamma$, 
      \label{item:pathcomponentsofIk}
   \item every element of $\Gamma$ contains at most 
     one path-connected component of $I_k$, and \label{item:onlyonecomponent}
   \item \label{item:complicated}
     if $c\subset I$ is a simple closed
     curve which is completely contained in some element
     of $\Gamma$, then either $c\subset I_k$ or $c\cap I_k = \emptyset$.
  \end{enumerate}

  Then $\Gamma$ is the set of path-connected components of $I$.
 \end{prop}
 \begin{remark}
  In the cases of interest to us, no element of $\Gamma$ will contain
   a simple closed curve, so (\ref{item:complicated}) is automatically
   satisfied.
 \end{remark}
 \begin{proof} 
  Since $I$ is Hausdorff, path-connectivity is the same as
   arc-connectivity. Let $\alpha:[0,1]\to I$ be any arc; we must show that
   $I$ is contained in a single element of $\Gamma$.

  For every $t\in [0,1]$, the point $\alpha(t)$ belongs to some
   $\gamma=\gamma_t\in \Gamma$; let us denote by
   $J(t)$ the component of $\alpha^{-1}(\gamma_t)$ containing $t$.  
   For every $k$, we also set $J'_k := \alpha^{-1}(I_k)$ (note that
   $J'_k$ is compact), and
   \[ J_k := \bigcup_{t\in J'_k} J(t). \]
   Our goal is to show that the hypotheses of Proposition
   \ref{prop:sigmaconnectedness} apply to $K=[0,1]$ and $A_k=J_k$.
   To do so, we need to check a number of simple
   properties.

 \begin{claim}[Claim 1]
    For every $t\in[0,1]$, the set $J(t)$ is compact.
 \end{claim}
 \begin{subproof}
  $J(t)$ is an interval; let
   $a$ and $b$ be the endpoints of $J(t)$. Then 
   $\gamma_t \cup \gamma_a$ is path-connected, since it contains
   the curve $\alpha\bigl([a,b)\bigr)$ containing $\alpha(a)$
   and $\alpha(t)$. 
   By the assumption on $\Gamma$, this implies $\gamma_t=\gamma_a$,
   and thus $a\in J(t)$. Analogously,
   $b\in J(t)$.
 \end{subproof}

 \begin{claim}[Claim 2]
   Let $t\in J'_k$. Then $J(t)$ contains exactly one component of
   $J'_k$. 
 \end{claim}
 \begin{subproof}
   Otherwise, there is 
    some interval $(a,b)\subset J(t)\setminus J_k'$ with $a,b\in J_k'$.
    Since $\alpha(a), \alpha(b)\in \gamma_t$,
    it follows by
    (\ref{item:onlyonecomponent}) that they belong to the
    same path-connected component of $I_k$. Let $\beta\subset I_k$
    be an arc connecting $\alpha(a)$ and $\alpha(b)$. By 
    (\ref{item:pathcomponentsofIk}), 
    $\beta\subset\gamma_a$. 
    Thus 
    $\beta \cup \alpha\bigl((a,b)\bigr)\subset \gamma_a$ 
    is a simple closed curve
    which intersects $I_k$ but is not contained in $I_k$, contradicting
    (\ref{item:complicated}). 
 \end{subproof}

 \begin{claim}[Claim 3]
   $J_k$ is compact, and, for every $t\in J_k$, $J(t)$ is a 
   connected component of $J_k$.
 \end{claim}
 \begin{subproof}
   Let $A$ be a component of $[0,1]\setminus J_k'\supset
   [0,1]\setminus J_k$. To
   fix ideas, let us suppose that neither $0$ nor $1$ is an endpoint of
   this interval, so that $A=(a,b)$ with $\cl{A}\cap J_k'=\{a,b\}$. 
   Then 
   $A\setminus J_k = A \setminus (J(a)\cup J(b))$.

  Since $J(a) \cup J(b)$ is compact by Claim 1, $A\setminus J_k$ is open. 
   Writing $[0,1]\setminus J_k = \bigcup A\setminus J_k$, where
   the union is taken over all intervals $A$ as above, it follows that
   $J_k$ is compact. 
   Furthermore, by Claim 2, $b\notin J(a)$. Thus $J(a)\neq J(b)$, so
   $A\setminus J_k \neq \emptyset$.
   Thus any two components of 
   $J_k'$ are separated by a point of $[0,1]\setminus J_k$. 

   Now let $t\in J_k$, and consider some point
   $t_1\in J_k\setminus J(t)$. Then
   $J(t)$ and $J(t_1)$ each contain a connected component of
   $J_k'$. These are separated by a point of $[0,1]\setminus J_k$, which
   consequently also separates $J(t)$ and $J(t_1)$. Thus $t$ and $t_1$
   belong to different connected components of $J_k$. So
   $J(t)$ is the connected component of $J_k$ containing $t$, completing
   the proof of Claim 3.
 \end{subproof}


 Claim 3 implies that every connected component of $J_k$ is also
  a connected component of $J_{k+1}$. Thus we can apply Proposition
  \ref{prop:sigmaconnectedness} to see that
  $J_k = [0,1]$ for some $k$. (Note that we can
  assume without loss of generality that $J_1\neq\emptyset$.)
  By Claim 3, this means that
   $\alpha\bigl([0,1]\bigr) = \alpha(J(0)) \subset \gamma_0$,
 as required.
 \end{proof}

 \begin{cor} \label{cor:dynamicrays}
  Let $\kappa\in\C$. Then every path connected component of
  $I(\Ek)$ is
  \begin{enumerate}
   \item a dynamic ray at a slow address or \label{item:slowray}
   \item a dynamic ray at a fast address together with its
          escaping endpoint, or \label{item:fastray}
    \item (if $\kappa\in I(\Ek)$) an iterated preimage component of
      the ray (possibly with endpoint) containing the
         singular value. \label{item:brokenray}
  \end{enumerate}
 \end{cor}
 \begin{remark}[Remark 1]
  According to the remark after Proposition \ref{prop:classification},
   every component $\gamma$
   of the last type is a curve to $\infty$. More precisely, according to
   whether the address $\s^0$ associated to the singular orbit is
   fast, periodic or slow but not periodic, $\gamma$ is 
   an arc with one finite escaping endpoint, an arc tending to infinity
   in both directions, or an injective curve to infinity (which may or
   may not land at a point $z\notin I(\Ek)$), respectively. 
 \end{remark}
 \begin{remark}[Remark 2]
  Note that, while different rays form different path-connected components of
   $I(\Ek)$, they may belong to the same path-connected component of
   the \emph{Julia set} $J(\Ek)$ (which is the closure of $I(\Ek)$). 
   Indeed, often $J(\Ek)$ will be the entire complex plane, and
   even otherwise several rays 
   may have a common (non-escaping) landing point.

  Similarly, different rays may belong to the same
   \emph{connected} (rather than path connected)
   component of $I(\Ek)$, e.g.\ in the postsingularly finite case
   \cite{misindecomposable}.
 \end{remark}
 \begin{proof} 
  Let $\Gamma$ be the set of all curves of the
  types (\ref{item:slowray}) to
  (\ref{item:brokenray}) listed above, which (by Proposition
  \ref{prop:classification}) exhausts all of $I(\Ek)$.
  Then the
   union of two different elements of $\Gamma$ is never path-connected,
  and no element of $\Gamma$ contains a simple closed curve.
  (This could only fail if there was some slow address
   $\s\in\Sequ_0$ for which $\gsk$ lands at an escaping point; this
   is impossible by \cite[Corollary 6.9]{expescaping}.)

  In \cite[Theorems 4.2 and 4.3]{topescapingnew}, 
   it was shown that there exists a relatively closed
   subset 
   $Y\subset I(\Ek)$ with the following properties.
   \begin{enumerate}
    \item $\Ek(Y)\subset Y$,
    \item for all $z\in I(\Ek)$, there is an $n\geq 0$ with
            $\Ek^n(z)\in Y$, and
    \item every connected component $C$ of $Y$ is the ``tail'' of some
      dynamic ray
      $\gs$, in the following sense:  either
      $C=\gs(X_{\s})$ or $C=\gs\bigl([t,\infty)\bigr)$ 
      for some $t>\ts$. 
      \label{item:straightbrush}
   \end{enumerate}
  (In fact, the set $Y$ is a ``straight brush''
   in the sense of Aarts and Oversteegen \cite{aartsoversteegen}.)

  It follows easily from (\ref{item:straightbrush})
   that every path-connected component of
   $I_n := \Ek^{-n}(Y)$ is contained in a single element of
   $\Gamma$, and that the intersection of
   $I_n$ with any element of $\Gamma$ is path-connected. 
   Thus we can apply Proposition
   \ref{prop:pathcomponentsarecurves} to $\Gamma$ and the sequence
   $(I_n)$, which completes the proof. 
 \end{proof}

 \begin{rmk}[Alternative definition of dynamic rays]
 By Corollary \ref{cor:dynamicrays},
  the following simple definition of dynamic rays 
  is equivalent to the original one:

  Let $\Ek$ be an exponential map and suppose that the singular value $\kappa$
   does not escape. A \emph{dynamic ray} of $\Ek$ is a maximal injective
   curve $\gamma\colon(0,\infty)\to\C$ in the escaping set%
   \footnote{As defined here, the parametrization 
   of $\gamma$ need not be related to potentials as used earlier; we include the parametrization
   over $(0,\infty)$ only in order to exclude the endpoint, if any, from the definition of dynamic ray.}.
 \end{rmk}

\section{Accumulation Behavior of Parameter Rays}
\label{sec:landingbehavior}

 In the proof of Corollary \ref{cor:dynamicrays}, we used the fact that
  a dynamic ray at a slow address cannot land at an escaping point. 
  To complete the proof of
  Theorem \ref{thm:main}, 
  we shall show the analogous statement in the parameter plane, 
  which requires two more facts about parameter rays. The first of these 
  is a result from
  \cite{expcombinatorics}. There the structure of parameter space
  was exploited to show that any two parameter rays 
  which accumulate at a common point must share the same
  ``kneading sequence'' (a combinatorial invariant), and
  in particular satisfy the following combinatorial restriction.

 \begin{prop}[\protect{\cite[Corollary A.2]{expcombinatorics}}]
  \label{prop:rayslandingtogether}
  Suppose that $\s^1,\s^2\in\Sequ_0$ such that
   $\PRS{\s^1}$ and $\PRS{\s^2}$ have a common accumulation point.
   Then $|s^1_k - s^2_k|\leq 1$ for all $k\geq 1$. \qedd
 \end{prop}

 Our second ingredient shows that a parameter ray cannot
  land at an escaping parameter which is not the endpoint of a parameter ray.
  (This is a special case of some
  results obtained in \cite{topescapingnew} regarding the 
  continuous dependence of 
  parameter rays on their external addresses: every ray is locally
  uniformly approximated by other parameter rays, so that the endpoint
  is the only possible accessible point.)

 \begin{prop}[{\cite[Corollary 11.3]{topescapingnew}}]
  \label{prop:raycont}
  Let $\s\in\Sequ_0$, and let $L_{\s}$ denote the
   accumulation set of $\PR$. 
   If 
    $\PRS{\addst}(t)\in L_{\s}$ for some $\addst\in\Sequ_0$ and
    $t>t_{\addst}$, then 
    $\PRS{\addst}\bigl((t_{\addst},t)\bigr)\subset
     L_{\s}$. \qedd
 \end{prop}

 \begin{cor} \label{cor:slowraysdontland}
  Suppose that $\s$ is a slow address. Then $\PR$ does not land at
  an escaping parameter.
 \end{cor}
 \begin{proof}
 By Proposition \ref{prop:raycont}, $\PR$ cannot land
  at a parameter of the form $\PRS{\addst}(t)$ with 
  $t>t_{\addst}$. So suppose, by contradiction, that
  $\PR$ lands at $\PRS{\addst}(t_{\addst})$ for some fast external
  address $\addst$. 
  Applying Proposition \ref{prop:rayslandingtogether},
  we see that
  $|\st_n - s_n| \leq 1$ for all $n$. This is impossible since
  $\s$ is slow and $\addst$ is fast.
 \end{proof}

 \begin{thm} \label{thm:pathcomponents}
  The path-connected components of $\I$ are exactly the sets
  $G_{\s}(X_{\s})$ with $\s\in\Sequ_0$.
 \end{thm}
 \begin{remark}
  This completes the proof of Theorem \ref{thm:main}.
 \end{remark}
 \begin{proof}
  Set $\Gamma := \{\PR(X_{\s}):\s\in\Sequ_0\}$; then
   the elements of $\Gamma$ are pairwise disjoint. Again, 
   the union of two such elements is 
   never path-connected, and no parameter ray
   contains a simple closed curve, by
   Corollary \ref{cor:slowraysdontland}. 

  We again wish to
   apply Proposition \ref{prop:pathcomponentsarecurves} to see that
   $\Gamma$ is the set of path-connected components of $\I$. To do so,
   we will construct
   an exhaustion of $\I$ by a sequence of
   ``straight brushes'', similarly as in the dynamical plane. (Compare
   \cite[Theorem 11.5]{topescapingnew}.) This is somewhat more
   complicated in parameter space, so let us begin by outlining
   the strategy.
   We define closed sets $\I'_n$, each of which is homeomorphic to a subset
   of a straight brush. If $\I_n$ is the union of all unbounded components 
   of $\I'_n$, then $\I_n$ is a straight brush, each of whose components
   is contained in exactly one element of $\Gamma$. To complete the proof,
   it remains to
   show that the sets $\I_n$ also exhaust $\I$, which is
   essentially a compactness argument. The details
   follow. 

  We begin by defining, for each $n\in\N$,
  \[ {\I'_n} := \{\kappa\in\I:
               \re \Ek^j(\kappa)\geq K(\kappa)
               \text{ for all $j\geq n$}\}, \]
  where $K(\kappa) := \log^+(|\kappa|) + 10$. 
  Clearly each ${\I_n}'$ is a closed subset of $\I$; furthermore
   ${\I}'_{n}\subset {\I'_{n+1}}$ and
   $\I = \bigcup {\I}'_n$.

   Let us endow the space 
   $\widehat{X} := X\cup\{\infty\}$ 
   with the ``straight brush''
   topology from \cite{topescapingnew}
   (compare the remark after Proposition \ref{prop:classification}). 
   As in the `ray-wise discrete' topology we have considered so far,
   every connected component of $X$ 
   is of the form
   $\{\s\}\times X_{\s}$, and every
   closed connected subset of $\widehat{X}$ is path-connected.

   Set $\Ih:= \I\cup\{\infty\}$, and
   consider the map 
   $\phi:\Ih \to \widehat{X}$ which takes $\infty$ to $\infty$ and
   maps every parameter
   $\kappa\in\I$ to the unique pair $(\s,t)\in X$ with
   $\kappa = \gsk(t)$. It follows from
   \cite[Corollary 4.4]{topescapingnew} that 
   the map $\phi$ is continuous on ${\I'_n}\cup\{\infty\}$
   for every $n$.%
      \footnote{With the notation after Corollary \ref{cor:classification}, 
      $\phi = \PRconj^{-1}$. Using the bijectivity of
      $\phi$, one can show that
      $\phi$ is actually a homeomorphism when restricted to
      $\I'_n$; however, both the definition
      of $\phi$ and its continuity on the mentioned
      sets are independent of Corollary \ref{cor:classification}.}%
  \ So
   any connected set $C\subset \I'_n$ is mapped
   into a connected component of $X$ under $\phi$, and 
   is thus a piece of
   a single parameter ray. If $C$ is unbounded, then 
   $C$ is a \emph{tail} of this ray 
   (in the same sense as in the proof of Corollary \ref{cor:dynamicrays}).
 
  Let $\Ih_n$ denote the component of
   ${\I'_n}\cup\{\infty\}$ containing $\infty$; by 
   the property of the
   topology of $\widehat{X}$ mentioned above, the set
   $\Ih_n$ is path-connected.
   We set
   $\I_n := \Ih_n\setminus\{\infty\}$.
   Then $\I_n$ 
   is a closed subset of $\I$. Since $\Ih_n$ is path-connected,
   every connected component of $\I_n$ is unbounded, and is thus
   a tail of some parameter ray.

 We can now apply Proposition \ref{prop:pathcomponentsarecurves} to 
  see that the path-connected
  components of $\widetilde{\I} := \bigcup \I_n$ are exactly the
  sets $\gamma\cap \widetilde{\I}$ with $\gamma\in\Gamma$. 
  It remains to show that $\widetilde{\I} = \I$; i.e., that
    $G_{\s}(t) \in \widetilde{\I}$
    for all $\s\in\Sequ_0$ and $t\in X_{\s}$.
  So let $\s\in\Sequ_0$ and consider the (closed) sets
  \begin{align*}
    A_n &:= \{t\in X_{\s}: G_{\s}(t)\in {\I'_n}\} \quad \text{ and} \\
    B_n &:= \{t\in X_{\s}: G_{\s}(t)\in \I_n\}.
  \end{align*}
  By \cite[Theorem 3.2]{markusdierk}, the set $B_1$ is nonempty. Note 
   that $B_n$
   is  the unbounded component of $A_n$ for all $n$, and that 
   $X_{\s} = \bigcup A_n$. We need to show that also
   $X_{\s} = \bigcup B_n$, which will be achieved through
   a standard connectivity
   argument. 

 \begin{claim}
  For every $T\in X_{\s}$, there exists $n\in\N$ such that
   $A_n$ is
   a neighborhood of $T$ in $X_{\s}$. 
 \end{claim}
 \begin{subproof}
   It follows easily from \cite[Theorem 4.2]{topescapingnew} that we
   can find a natural number $n=n(T)$ with the following property:
   if $\kappa$ is sufficiently close to $G_{\s}(T)$ and 
    $t\in X_{\s}$ is sufficiently close to $T$, then 
   \[ \re \Ek^m(g^{\kappa}_{\s}(t)) > K(G_{\s}(T)) + 1 \]
   for all $m\geq n$. 
  Also, $K(\kappa)$ is continuous in $\kappa$, and 
   $G_{\s}(t)$ is continuous in $t$. So if $t$ is close enough to $T$,
   then $\kappa := G_{\s}(t)$ satisfies
   \[ \re \Ek^m(\kappa) = \re \Ek^m(g^{\kappa}_{\s}(t)) >
       K(G_{\s}(T)) + 1 \geq K(\kappa) \]
   for all $m\geq n$. Therefore
   $t \in A_n$, completing the proof of the Claim.
 \end{subproof}


 For every $t\in X_{\s}$, define
  \[ B(t) := \bigcup\{C: C\text{ is the connected component of
                         $A_n$ containing $t$ for some $n\in\N$}\}. \]
 By the Claim, $B(t)$ is open for every $t$, so the
  $B(t)$ form a partition of $X_{\s}$ into open sets. 
  Since $X_{\s}$ is connected, this means that
  $B(t)=X_{\s}$ for all $t$. On the other hand, for $t_0\in B_1$ and
  $n\in\N$, 
  $B_n$ is the connected component of $A_n$ containing $t_0$. Thus
  $\bigcup B_n = B(t_0) = X_{\s}$.
  This completes the proof.
 \end{proof}

\section*{Acknowledgments}
We would like to thank the referee for useful suggestions, and 
  Vlad Vicol for an interesting historical discussion.


\nocite{jackdynamicsthird}
\nocite{alexmisha}
\nocite{dghnew1,dghnew2}
\bibliographystyle{hamsplain}
\bibliography{../../Biblio/biblio}

\end{document}